\documentclass[final,1p,times]{elsarticle}
\usepackage{amssymb,amsmath,amsthm}
\usepackage[all]{xy}
\usepackage{hyperref}

\usepackage{amscd}

\usepackage{enumerate}

\usepackage{graphicx}

\newcommand{\bx}{{\mathbf x}}

\newcommand{\C}{{\mathbb C}}

\DeclareMathOperator{\sech}{sech}
\DeclareMathOperator{\cond}{cond}

\makeatletter
\def\Ddots{\mathinner{\mkern1mu\raise\p@
\vbox{\kern7\p@\hbox{.}}\mkern2mu
\raise4\p@\hbox{.}\mkern2mu\raise7\p@\hbox{.}\mkern1mu}}
\makeatother

\newtheorem{condthm}{Theorem}
\newtheorem{definition}{Definition}
\journal{ }

\begin{document}
\begin{frontmatter}
\title{Norm-preserving discretization of integral equations
for elliptic PDEs with internal layers I: the one-dimensional case}
\author{Travis Askham and Leslie Greengard}
\address{Courant Institute,
New York University, New York, NY 10012.}

\begin{abstract}
We investigate the behavior of 
integral formulations of variable coefficient elliptic
partial differential equations (PDEs) in the presence of steep internal layers.
In one dimension, the equations that arise can be solved 
analytically and the condition numbers estimated in 
various $L^p$ norms. 
We show that high-order accurate Nystr\"{o}m discretization leads to 
well-conditioned 
finite-dimensional linear systems if and only if 
the discretization is both norm-preserving in a correctly 
chosen $L^p$ space {\em and} adaptively refined in the internal layer.
\end{abstract}
\end{frontmatter}


\section{Introduction}

A number of problems in computational physics require the solution of
divergence-form elliptic equations
\begin{equation}
\label{pde1}
 \nabla \cdot (\epsilon(\bx) \nabla u(\bx)) = f(\bx)  
\end{equation}
where $\epsilon(\bx)$ is a scalar function with steep internal layers
in a domain $\Omega$. We assume for the sake of concreteness 
that $u(\bx)$ satisfies a Dirichlet boundary condition
\begin{equation}
\label{pdebc}
u(\bx) = g(\bx)  
\end{equation}
for $\bx \in \partial\Omega$,
but the basic approach outlined below applies equally well
to other types of boundary conditions.
Equations of the form (\ref{pde1}) arise,
for example, in fluid dynamics \cite{bell,lions}, 
where $\epsilon(\bx)$ is the inverse of the fluid density and in 
semiconductor device simulation \cite{Markowich}, 
where $\epsilon(\bx)$ can be either the semiconductor permittivity,
or a complicated function determined by electron and hole mobilities
and diffusion coefficients.
They also arise in phase field models for microstructure evolution
in materials science \cite{phasefield}.
When $\epsilon$ is piecewise constant, boundary
integral equation methods are well-known to be extremely effective
(see, for example, \cite{glee,gmoura,helsing2,martinsson,white}).
When $\epsilon$ is smooth but has a steep internal layer,
however, the domain itself must be discretized. 
In that setting, it is most common to use 
finite difference or finite element approximations
based on the partial differential equation itself
\cite{brenner,leveque,sptext}.

Volume integral equations can also been used for problems
such as (\ref{pde1}). There is a substantial literature in this area, 
which we do not attempt to review, except to observe that there are a 
variety of analytic methods which can be used to derive 
integral formulations, a variety of numerical methods  which can be used
for their discretization, and a variety of fast algorithms which can be used
for iterative or direct solution
\cite{chandrasekaran,CHEN,Chewbook,CK2,helsing,gillman,hackbusch,ho,kottmann,vrvol}.

In this paper, we focus on the behavior of volume integral methods
in one dimension, where the divergence form equation reduces to
\begin{equation}
\label{ode1}
\frac{\partial}{\partial x} 
\left( \epsilon(x) \frac{\partial u}{\partial x} \right) = f \, .
\end{equation}

For the sake of simplicity, we assume
the solution is subject to homogeneous Dirichlet conditions on the 
interval $[a,b]$, that is $u(a) = u(b) = 0$.
We assume that $\epsilon(x)$ is positive, smooth and bounded,
but may have steep gradients, so that its derivative $\epsilon_x(x)$
can be arbitrarily large, corresponding to an
internal layer. Without care, this can lead to arbitrarily 
badly conditioned linear
systems. 
While there is some literature on analyzing the
conditioning of second kind integral
equations (see, for example, \cite{ahues,ikebe}), the 
influence of the choice of $L^p$ space has received 
relatively little attention. 
Here, we show that a combination of adaptivity and 
a suitable {\em norm-preserving} discretization, to be defined below, 
leads to condition numbers that depend only weakly on 
$\epsilon_x$. In particular, we show that for a Lippmann-Schwinger
type integral equation with the second
derivative $u_{xx}$ as the unknown, a discretization that 
is norm-perserving in $L^1$ leads to nearly optimal schemes. 

Our work was motivated, in part, by Bremer's analysis of boundary
integral equations for scattering problems in the 
presence of corners \cite{bremer12}. He showed that naive 
Nystr\"{o}m discretization leads to ill-conditioned linear systems,
but that suitable $L^2$-weighting corrects the difficulty both in theory
and in practice.

\section{The integral equation}

There are several standard methods for converting the ordinary 
differential equation (\ref{ode1}) to an integral equation, 
typically making use of the Green's function $G(x,t)$ that satisfies
\[ \frac{d^2}{dx^2} G(x,t) = \delta(x-t),\qquad
G(a,t) = G(b,t) = 0 \, . \]
It is well-known \cite{guenther} and easy to verify that 
\begin{equation}
\label{gdef}
G(x,t) = \left \{ \begin{array}{rcl} (x-a)(t-b)/(b-a) & \mbox{if} & x < t \\
                       (x-b)(t-a)/(b-a) & \mbox{if} & x \geq t
                      \end{array} \, . \right.
\end{equation}
Rewriting the equation (\ref{ode1}) in the form
\begin{equation}  
\label{ode2}
u_{xx} + \frac{\epsilon_x}{\epsilon} u_x =  \frac{f}{\epsilon}
\end{equation}
and representing the solution as 
\begin{equation}
u(x) = \int_a^b G(x,t) \sigma(t) \, dt \, ,
\label{urepsigma}
\end{equation}
we obtain an integral equation
for the unknown density $\sigma$:
\begin{equation}
\label{inteq1}
 \sigma(x) + 
\frac{\epsilon_x(x)}{\epsilon(x)} \int_a^b G_x(x,t) \sigma(t) \, dt = g(x) \, ,
\end{equation}
or
\begin{equation}
\label{inteq1a}
(I + K_1) \sigma(x) = g(x) 
\end{equation}
where $g = f/\epsilon$ and 
\[  K_1 \sigma(x) =
\frac{\epsilon_x(x)}{\epsilon(x)} \int_a^b G_x(x,t) \sigma(t) \, dt \, . 
\]
Alternatively, one can rewrite (\ref{ode1}) in the form
\begin{equation}  
\label{ode3}
(\epsilon u)_{xx} - (\epsilon_x u)_x =  f \, .
\end{equation}
Integrating (\ref{ode3}) against $G(x,t)$ yields
\begin{equation}
\label{inteq2}
 u(x) + \frac{1}{\epsilon(x)} \int_0^1 G_x(x,t)
(\epsilon_x(t) u(t)) \, dt = 
\frac{1}{\epsilon(x)} \int_0^1 G(x,t) f(t) dt 
\end{equation}
or
\begin{equation}
\label{inteq2a}
 (I + K_2)u(x) =
\frac{1}{\epsilon(x)} \int_a^b G(x,t) f(t) dt \, ,
\end{equation}
where 
\[  K_2 u(x) = \frac{1}{\epsilon(x)} \int_0^1 G_x(x,t)
(\epsilon_x(t) u(t)) \, dt \, .
\]
The principal difference between (\ref{inteq1}) and 
(\ref{inteq2}) is that, in the former, 
$\sigma(x) = u_{xx}(x)$ is the unknown
while, in the latter, $u(x)$ is the unknown.
Both are Fredholm equations of the second kind.

\subsection{Analytic solution of the integral equation}

For the sake of simplicty, let us assume in this section that $[a,b]=[0,1]$.
>From the original ODE, we have  
\begin{align}
   (\epsilon(x) u_x(x))_x &= g(x) \epsilon(x) \nonumber \\
   \epsilon(x) u_x(x) &= \int_0^x g(t) \epsilon(t) \, dt + \epsilon(0) u_x(0) 
\nonumber \\
   u_x(x) &= \frac{1}{\epsilon(x)} \int_0^x g(t) \epsilon(t) \, dt + \frac{\epsilon(0)u_x(0)}{\epsilon(x)} 
\label{uxeq}
  \end{align}
Using the fact that $\sigma = u_{xx}$, we may write
\begin{equation}
   \sigma(x) = g(x) - \frac{\epsilon_x(x)}{\epsilon(x)^2} \left ( \int_0^x g(t) \epsilon(t) \, dt + \epsilon(0)u_x(0) \right ) \, .
\label{resolve1}
\end{equation}
To remove the $\epsilon(0)u_x(0)$ term from the expression, we integrate
the equation (\ref{uxeq}).
\[
   u(1)-u(0) = \int_0^1 \frac{1}{\epsilon(x)} \int_0^x g(t) \epsilon(t) \, dt \, dx + \epsilon(0)u_x(0) \int_0^1 \frac{1}{\epsilon(x)} \, dx 
\]
so that
\begin{equation}
\label{eps0ux0}
   \epsilon(0)u_x(0) = - \frac{\int_0^1 \frac{1}{\epsilon(x)} \int_0^x g(t) \epsilon(t) \, dt \, dx} { \int_0^1 \frac{1}{\epsilon(x)} \, dx} \, .
\end{equation}
Letting $A_1 = I + K_1$ denote the operator applied to $\sigma$ on 
the left-hand side of (\ref{inteq1}), we now have an expression for
its inverse in the form $A_1^{-1} = I - R_1$.
>From (\ref{resolve1}) and (\ref{eps0ux0}),  
\[ \sigma(x) = g(x) - \frac{\epsilon_x(x)}{\epsilon(x)^2} \left ( \int_0^x g(t) \epsilon(t) \, dt  - \frac{\int_0^1 \frac{1}{\epsilon(s)} \int_0^s g(t) \epsilon(t) \, dt \, ds} { \int_0^1 \frac{1}{\epsilon(s)} \, ds} \right ) \, .
\]
>From this, it is straightforward to obtain the following formula for the 
resolvent kernel $R_1$:
\begin{equation} \label{resolvent1} R_1(x,t) =  \frac{\epsilon_x(x)}{\epsilon(x)^2} \left(
H(x-t) \epsilon(t) - \frac{\epsilon(t)}{\int_0^1 \frac{1}{\epsilon(s)} \, ds}
  \int_t^1 \frac{1}{\epsilon(s)} ds
\right) \, ,
\end{equation}
where $H(x)$ is the standard Heavyside function.

Letting $A_2 = I + K_2$ denote the operator applied to $u$ on 
the left-hand side of (\ref{inteq2}), a similar calculation
yields an expression for
its inverse in the form $A_2^{-1} = I - R_2$.
In this case, $R_2$ is 
\begin{equation} \label{resolvent2} R_2(x,t) =  -\frac{\epsilon_x(t)}{\epsilon(t)^2} \left(
H(x-t) \epsilon(t) - \frac{\epsilon(t)}{\int_0^1 \frac{1}{\epsilon(s)} \, ds}
  \int_0^x \frac{1}{\epsilon(s)} ds
\right) \, .
\end{equation}

Having analytic expressions for the resolvent kernels permits us to
obtain simple estimates for the condition number of the 
operators $A_1$ and $A_2$ acting on $L_p$ spaces for $1 \leq p \leq \infty$.
It is worth noting an 
important difference between the two resolvent kernels:
the term $\epsilon_x/\epsilon^2$ in (\ref{resolvent2}) is evaluated 
at $t$ rather than $x$. It is integrated when applying the 
inverse operator:
\begin{equation}
 \label{sol2} 
u(x) = h(x) + \int_0^x h(t) \frac{\epsilon_x(t)}{\epsilon(t)} \, dt - \frac{ \int_0^1 h(s) \frac{ \epsilon_x(s)}{\epsilon(s)} \, ds}{\int_0^1 \frac{1}{\epsilon(s)} \, ds} \int_0^x \frac{1}{\epsilon(t)} \, dt 
\end{equation}

\section{Integral Equation Operator Bounds}
We wish to characterize functions $\epsilon(x)$ that (a) are fairly flat on some subinterval of $[a,b]$ and 
(b) are uniformly bounded from above and below.
These conditions are formalized as follows:

\begin{definition} \label{propdef}
Let $\mathcal{E}$ denote a family of functions on the interval $[a,b]$.
\begin{itemize}
\item $\mathcal{E}$ satisfies {\em Property 1} if there 
exists $0 \leq \delta \ll 1$ and a constant $c > 0$ such that, for each 
$\epsilon \in \mathcal{E}$, there is a neighborhood 
$V = B( \zeta(\epsilon),c) \subset [a,b]$ such that 
\[ \left \| \frac{\epsilon_x}{\epsilon} \cdot 1_V \right \|_p \leq \delta \left \| \frac{\epsilon_x}{\epsilon} \right \|_p \]
for all $1 \leq p \leq \infty$.
\item $\mathcal{E}$ satisfies {\em Property 2} 
if $m > 0$ and $M<\infty$ where  
\[ m = \inf_{\epsilon \in \mathcal{E}} \left[ \min_{x \in [a,b]} 
\epsilon(x) \right] \mbox{ and } 
M = \sup_{\epsilon \in \mathcal{E}} \left[ \max_{x \in [a,b]} 
\epsilon(x) \right] 
\]
\end{itemize}
\end{definition}

We then have the following result on the condition number of the 
operator $A_1$, the Fredholm operator on
the left-hand side of (\ref{inteq1}).
   
\begin{condthm} \label{condthm1}
Let $\mathcal{E}$ be a family of functions satisfying Properties 1 and 2. 
Then
\[ C_1 \left( \left\| \frac{ \epsilon_x}{\epsilon }  \right\|_p^2 \right) 
-1 \leq \mbox{cond}_p(A_1(\epsilon)) \leq C_2 
\left( \left\| \frac{ \epsilon_x}{\epsilon }  \right\|_p^2 \right) + 1 \, ,
\]
where $\mbox{cond}_p(A_1(\epsilon))$ is the condition number of 
$A_1(\epsilon)$ as an operator from $L^p[a,b] \to L^p[a,b]$ for 
$1 \leq p < \infty$ and as an operator from $L^\infty[a,b] \cap C[a,b] 
\to L^\infty[a,b] \cap C[a,b]$ for $p = \infty$.
\end{condthm}

A proof can be found in the Appendix. Theorem \ref{condthm1} 
gives us a sense of the qualitative behavior of $A_1(\epsilon)$ acting 
on $L^p$ spaces. In particular,
its condition number is well-controlled in $L^1$, even when
there are steep  internal layers (where $\epsilon_x/\epsilon$ can be large). 
In $L^1$, it is the total variation of $\epsilon$ that matters.
In the $L_\infty$ norm, on the other hand, the operator norm can be seen to 
be large by inspection.
A dual result can be obtained for the integral operator
$A_2(\epsilon) = I+K_2$ in (\ref{inteq2a}). 

\begin{condthm} \label{condthm2}
Let $\mathcal{E}$ be a family of functions satisfying Properties 1 and 2. 
Then
\[  C_1 \left( \left\| \frac{ \epsilon_x}{\epsilon }  \right\|_q^2 
\right) - 1 \leq \mbox{cond}_p(A_2(\epsilon)) \leq C_2
\left( \left\| \frac{ \epsilon_x}{\epsilon }  \right\|_q^2 \right) + 1 \, ,
\]
where $1/p+1/q=1$ and $\mbox{cond}_p(A_2(\epsilon))$ is the condition number 
of $A_2(\epsilon)$ as an operator from $L^p[a,b] \to L^p[a,b]$ for 
$1 \leq p < \infty$ and as an operator from $L^\infty[a,b] \cap C[a,b] 
\to L^\infty[a,b] \cap C[a,b]$ for $p = \infty$.
\end{condthm}
Since the condition number in $L^p$ depends on the $L^q$ norm of 
$\epsilon_x/\epsilon$ in this case, it is clear that the condition number of 
$A_2(\epsilon)$ will be modest in $L^\infty$ and very large in $L^1$
in the presence of internal layers.

\section{Norm-Preserving Discretization} \label{normpreserve}

In order to analyze the condition number of discretized
integral equations, it is convenient to introduce the following definition.

\begin{definition}
A mapping $\Phi:V \subset L^p[a,b] \to \C^n$ is said to be
{\em norm-preserving} if
\[ \|\Phi ( g) \|_{l^p} = \| g \|_{L^p[a,b]} \]
for all $g \in V$. 
\end{definition}
Let $A$ be an invertible, bounded integral operator mapping $V$ to $U$. 
We say that a matrix $A_h(V)$ is a norm-preserving discretization of 
$A$ on the subspace $V$ if there exist norm-preserving mappings 
$\Phi$ and $\Psi$ such that the diagram
\[ \begin{CD}
  V \subset L^p[a,b] @>A>> U \subset L^p[a,b] \\
  @VV\Psi V                @VV\Phi V \\
  \C^n @>A_h>> \C^n
  \end{CD} \]
commutes. 

In the Hibert space case $(p=2)$, it was shown in \cite{bremer12}
that inner product preserving 
discretizations have singular values which approximate those of the original 
operator. 
In the Banach space setting, it is easy to show something equally useful, 
namely that the condition number of a 
norm-preserving discretization approximates that of the original operator. 

For this, let $B|_W$ denote the restriction of an operator $B$ to a 
subspace $W$. Let $A$ be an invertible, bounded operator mapping $V$
to $U$, let $\Psi$, $\Phi$ be norm-preserving mappings and let 
$A_h$ be a norm-preserving discretization of $A$, as above.
Then,
\begin{align}   
\| A_h|_{\Psi(V)} \|_{l^p} &= \sup_{ v \in \Psi(V) } \frac{ \|A_h v \|_{l^p} }{\|v\|_{l^p} } = \sup_{ g \in V } \frac{ \| Ag \|_{L^p} }{ \|g\|_{L^p} } = \|A|_V \|_{L^p} \, ,  \label{opnorm}\\
\|A_h^{-1}|_{\Phi(U)} \|_{l^p} &=  \sup_{ w \in \Psi(V)} \frac{ \|w \|_{l^p} }{\|A_h w\|_{l^p} }  = \sup_{f \in V} \frac{ \|f\|_{L^p}}{\|Af\|_{L^p}} = \|A^{-1}|_U \|_{L^p}  \, .\label{inversenorm}
\end{align}
Thus, the condition number of $A_h$ restricted to $\Psi(V)$ and of
$A$ restricted to $V$ are the same.
  
\subsection{Norm-preserving Nystr\"{o}m discretizations}\label{nystrom}

We build (approximate) 
norm-preserving Nystr\"{o}m discretizations for $A$ by applying a 
quadrature rule
to the integral operator $A = I+K$:
\[ A f(x) = f(x) + \int_a^b K(x,y) f(y) \, dy \, . \]
For this, we assume that we are given an $n$-point quadrature rule
\[ \int_a^b f(x) \, dx \approx \sum_{k=1}^n f(x_k) w_k \, , \]
with positive weights. This induces a mapping 
$\Phi$: $L^p[a,b] \rightarrow \C^n$:
\begin{equation} 
\label{npmap} 
\Phi(f) = \left( \begin{array}{c} f(x_1) w_1^{1/p} \\ 
\vdots \\ f(x_n) w_n^{1/p} \end{array} \right) 
\end{equation}
If the quadrature rule is exact for functions of the form $|g|^p$ for 
$g \in V$ and $|f|^p$ for $f \in U$, then $\Phi$ is a norm-preserving 
mapping from $V$ into $\C^n$ and $U$ into $\C^n$. 
Further, suppose that the quadrature rule is exact for functions of the 
form $K(x,\cdot) g(\cdot)$ where $g \in V$, and that $A_h$ is given by
th Nystr\"{o}m discretization: 
\begin{equation} 
(A_h)_{ij} = \delta_{ij} + K(x_i,x_j) w_i^{1/p} w_j^{1-1/p} \, .
\label{ahijeq} 
\end{equation}
Then $A_h$ is norm-preserving, since
  \begin{align}   
[A_h \Phi(g)]_i &= g(x_i) w_i^{1/p} + w_i^{1/p} 
\sum_{j=1}^n K(x_i,x_j) w_j^{1-1/p} g(x_j) w_j^{1/p} \\
&= w_i^{1/p} \left( g(x_i) + \int_a^b K(x_i,y) g(y) \, dy \right) \, .
\end{align}
We note that discretization by {\em sampling}, i.e. where 
\[ \Phi(f) = \left( \begin{array}{c} f(x_1) \\ \vdots \\ 
f(x_n) \end{array} \right) \]
corresponds to 
a norm-preserving Nystr\"{o}m discretization on the space 
$L^\infty[a,b] \cap C[a,b]$. In particular,  suppose we 
let $V \subset L^\infty[a,b] \cap C[a,b]$ be equicontinuous and 
let $0<\delta \ll 1$. Then, by taking a fine enough mesh we can clearly
satisfy
\[ \|f\|_{L^\infty} = \|\Phi(f)\|_{l^\infty} (1+\delta) \]
for any $f \in V$. In short, the simplest
Nystr\"{o}m discretization, corresponding to sampling the unknown
on a grid, results in a discrete operator whose condition number 
approximates that of the continuous operator acting on 
$L^\infty[a,b] \cap C[a,b]$.  
    
\subsection{Discrete condition number estimates in alternate norms} 
\label{adaptivity}

Two aspects of norm-preserving discretations should be noted here.
First, the fact that a
discretized operator equation is well-conditioned in $l^p$ for some $p$ 
may not be very informative if we solve the finite-dimensional linear
algebra problem using a different norm.
Suppose, for example, that we wish to solve the equation 
(\ref{inteq1}), which is well-conditioned in $L^1$. After discretization
using (\ref{ahijeq}), it is well-conditioned in $l^1$ as well. 
However, if we use an iterative scheme such as GMRES \cite{GMRES},
we would like to ensure rapid convergence, which depends on the 
condition number in $l^2$.
(One could, of course, solve linear systems iteratively in $l^p$ spaces,
but the procedures are nonlinear and much more expensive.)

Fortunately, in finite dimensional spaces, norms and condition numbers
are all equivalent and satisfy simple relations \cite{GOLUB}. 
For instance,
\begin{equation}
\cond_2 (A_h) \leq n \cond_1(A_h) \, .
\end{equation}
Thus, if the system size is modest and we employ a norm-preserving 
discretization for $L^1$, we will have an acceptable bound on the $l^2$ 
condition number of the system matrix 
(\ref{ahijeq}).

A second, closely related, feature of norm-preserving discretizations is that
spatial adaptivity is {\em essential} for the choice of $L^p$ to have 
an impact.  One can see from (\ref{ahijeq}) that for 
a uniform mesh (with $w_i = h = \frac{1}{n}$ for all $i$),
the resulting matrix $A_h$ is the same for every $p$. 
Thus, if the continuous operator equation has a large condition 
number in $L^2$, the discretized equation will be ill-conditioned in $l^2$ 
as well. 

We will return to these issues in section \ref{sec:discussion}, following
an exploration of the behavior of the $l^1$, $l^2$ and $l^\infty$ 
discretizations on some model problems.

\section{Numerical Examples} \label{numerical}

To investigate the utility of the analysis outlined above, let us
first consider functions $\epsilon(x)$ in 
(\ref{ode1}) of the form
\begin{equation}
\label{epseq}
\epsilon_\delta (x) = 2+\tanh ( \delta ( x-x_0)) 
 \end{equation}
 on the interval $[0,2]$, where $x_0 \in (0,2)$. 
For large values of $\delta$, these functions have a steep internal layer 
centered at $x=x_0$. They are relatively flat away from the internal layer. and 
they are bounded in the range $[1,3]$. As a result,
the family 
\begin{equation}
  \label{epsfamily}
  \mathcal{E} = \left \{ \epsilon_\delta \in L^p : \delta \geq 10 \right \} 
 \end{equation}
satisfies Properties 1 and 2 as given in Definition \ref{propdef}. 
Note that the derivative 
$(\epsilon_\delta)_x = \delta \sech^2 ( \delta (x-x_0))$, so that
 
 \begin{align}
  \| (\epsilon_\delta)_x \|_p &= \left ( \int_0^2 \delta^p \sech^{2p}(\delta(x-x_0)) \, dx \right )^{1/p} \nonumber \\
  &\leq \delta \left ( \int_0^2 \sech^{2}(\delta(x-x_0)) \, dx \right )^{1/p} \nonumber \\
  &= \delta^{(1-1/p)} \left ( \tanh(\delta(2-x_0)) + \tanh( \delta x_0) \right )^{1/p} \label{epsxup} \\
  \| (\epsilon_\delta)_x \|_p &= \left ( \int_0^2 \delta^p \sech^{2p}(\delta(x-x_0)) \, dx \right )^{1/p} \nonumber \\
  &\geq 
  \frac{\delta}{2} \left (2 \cosh^{-1}(\sqrt{2})/ \delta \right)^{1/p} \nonumber \\
  &= C(p) \delta^{(1-1/p)} \label{epsxlow}
 \end{align}
Combining (\ref{epsxup}) with (\ref{epsxlow}) and the fact that 
the $\epsilon_\delta$ are uniformly bounded above and below, we have 
 
 \begin{equation}
  \left \| \frac{ ( \epsilon_\delta )_x}{\epsilon_\delta} \right \|_p = \Theta \left (\delta^{(1-1/p)} \right ) 
 \end{equation}
 for $1 \leq p < \infty$, using the standard ``Big Theta'' notation. 
It is straightforward to check that
\begin{equation}
  \left \| \frac{ ( \epsilon_\delta )_x}{\epsilon_\delta} \right \|_\infty = \Theta \left (\delta \right ).
\end{equation}
Letting $A_1(\epsilon)$ and $A_2(\epsilon)$ be the operators given 
by the left hand sides of (\ref{inteq1a}) and (\ref{inteq2a}), respectively,
and applying Theorem 1 to the family $\mathcal{E}$, we see that
\begin{align*}
\mbox{cond}_1(A_1(\epsilon_\delta)) &= \Theta(1), \\
\mbox{cond}_2(A_1(\epsilon_\delta)) &= \Theta(\delta), \\
\mbox{cond}_\infty (A_1(\epsilon_\delta)) &= \Theta(\delta^2).
\end{align*}
Likewise, we have  
\begin{align*}
\mbox{cond}_1(A_2(\epsilon_\delta)) &= \Theta(\delta^2), \\
\mbox{cond}_2(A_2(\epsilon_\delta)) &= \Theta(\delta), \\
\mbox{cond}_\infty (A_2(\epsilon_\delta)) &= \Theta(1).
\end{align*}
  
We discretize the integral equations (\ref{inteq1}) and (\ref{inteq2}), 
using a norm-preserving Nystr\"{o}m discretization scheme, as described in 
section \ref{nystrom}. 
For this, we adaptively refine the interval $[a,b]$ so that the function
$\epsilon(x)$ is well resolved with a piecewise Legendre polynomial
approximation to a user-specified precision.
More precisely, we use piecewise 16{\em th} order approximations,
and refine each interval until the quadrature error in integrating
$\epsilon$ is less than $10^{-15}$.
On each subinterval, we sample all functions involved $(u,\epsilon,f)$ at
the scaled Gauss-Legendre nodes of order $16$. 
We use the standard Gauss-Legendre quadrature weights scaled to each 
subinterval. 
Given these nodes and weights, the norm-preserving discretization 
(\ref{ahijeq}) in $L^p$ applied to equation (\ref{inteq1}) yields
\begin{equation}
 \sigma(x_i) w_i^{1/p} + \frac{\epsilon_x(x_i)}{\epsilon(x_i)} 
\sum_j G_x(x_i,x_j) w_j^{1-1/p} w_i^{1/p} 
\sigma (x_j) w_j^{1/p} = g(x_i) w_i^{1/p} \label{linsystem1} \, .
\end{equation}
Likewise, equation (\ref{inteq2}) yields 
\begin{equation}
 u(x_i) w_i^{1/p} + \frac{1}{\epsilon(x_i)} \sum_j G_x(x_i,x_j) 
\epsilon_x(x_j) w_j^{1-1/p} w_i^{1/p} u(x_i) w_i^{1/p} = h(x_i) 
w_i^{1/p} \label{linsystem2}
\end{equation}
where $h$ is simply the right hand side of (\ref{inteq2}). 
We will use $A_{1,p}(\epsilon)$ and $A_{2,p}(\epsilon)$ to denote 
the $p$-norm-preserving discretizations of these integral operators. 
Because the unknowns $\sigma$ and $u$ are weighted by $w_i^{1/p}$, 
we see that the entries of the discrete operators are given by

\begin{align}
 \left [ A_{1,p}(\epsilon) \right]_{ij} &= \delta_{ij} + \frac{\epsilon_x(x_i)}{\epsilon(x_i)} G_x(x_i,x_j) w_j^{1-1/p} w_i^{1/p} \nonumber \\
\left [ A_{2,p}(\epsilon) \right]_{ij} &= \delta_{ij} + \frac{\epsilon_x(x_j)}{\epsilon(x_i)} G_x(x_i,x_j) w_j^{1-1/p} w_i^{1/p} \nonumber
\end{align}

\subsection{Condition Numbers}

Using the family of functions $\mathcal{E}$ defined above, 
we may study the $l^p$ condition numbers of our discrete operators 
$A_{1,p}(\epsilon_\delta)$ and $A_{2,p}(\epsilon_\delta)$ for 
$p=1$, $2$, and $\infty$. 
Because of the norm-preserving discretization,
we expect $\mbox{cond}_1(A_{1,1}(\epsilon_\delta)) = \Theta(1)$,
$\mbox{cond}_2(A_{1,2}(\epsilon_\delta)) = \Theta(\delta)$, and
$\mbox{cond}_\infty (A_{1,\infty}(\epsilon_\delta)) = \Theta(\delta^2)$
since that is the behavior of the continous operators 
(Theorem \ref{condthm1}).
Similarly, we expect
$\mbox{cond}_1(A_{2,1}(\epsilon_\delta)) = \Theta(\delta^2)$, 
$\mbox{cond}_2(A_{2,2}(\epsilon_\delta)) = \Theta(\delta)$, and 
$\mbox{cond}_\infty (A_{2,\infty}(\epsilon_\delta)) = \Theta(1)$
(from Theorem \ref{condthm2}).

\begin{figure}[ht]
 \centering
    \includegraphics[width=.3\textwidth]{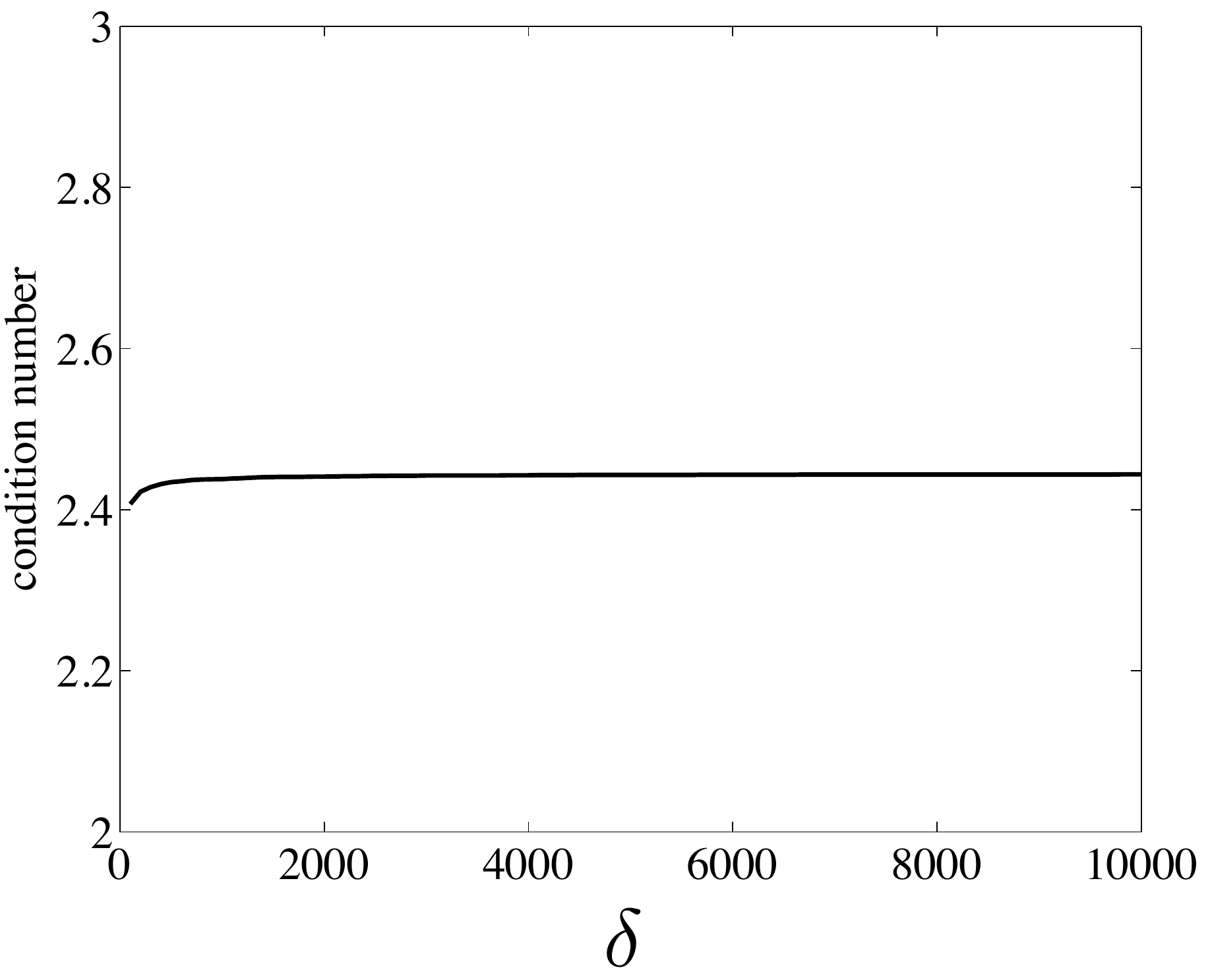}
    \includegraphics[width=.31\textwidth]{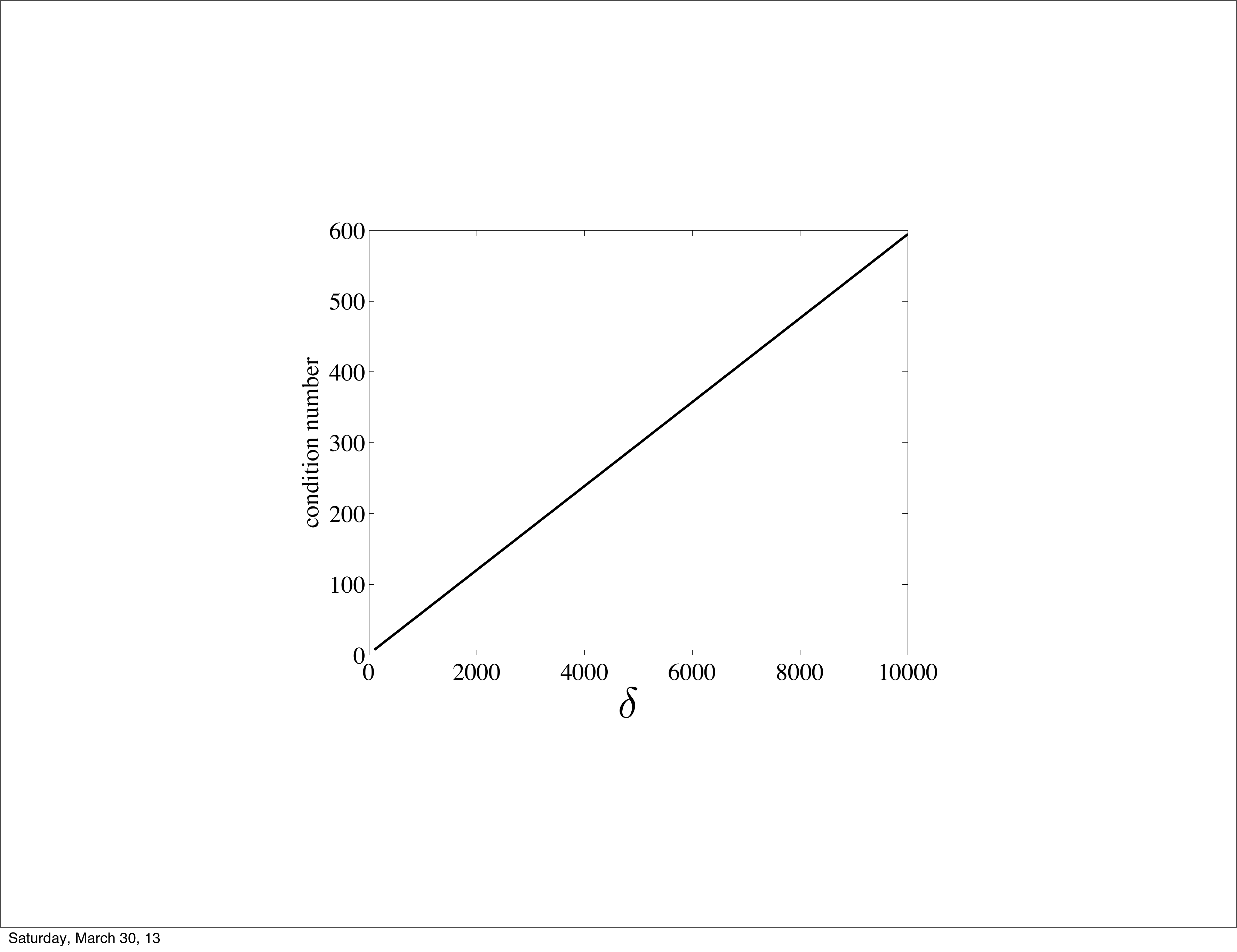}
    \includegraphics[width=.333\textwidth]{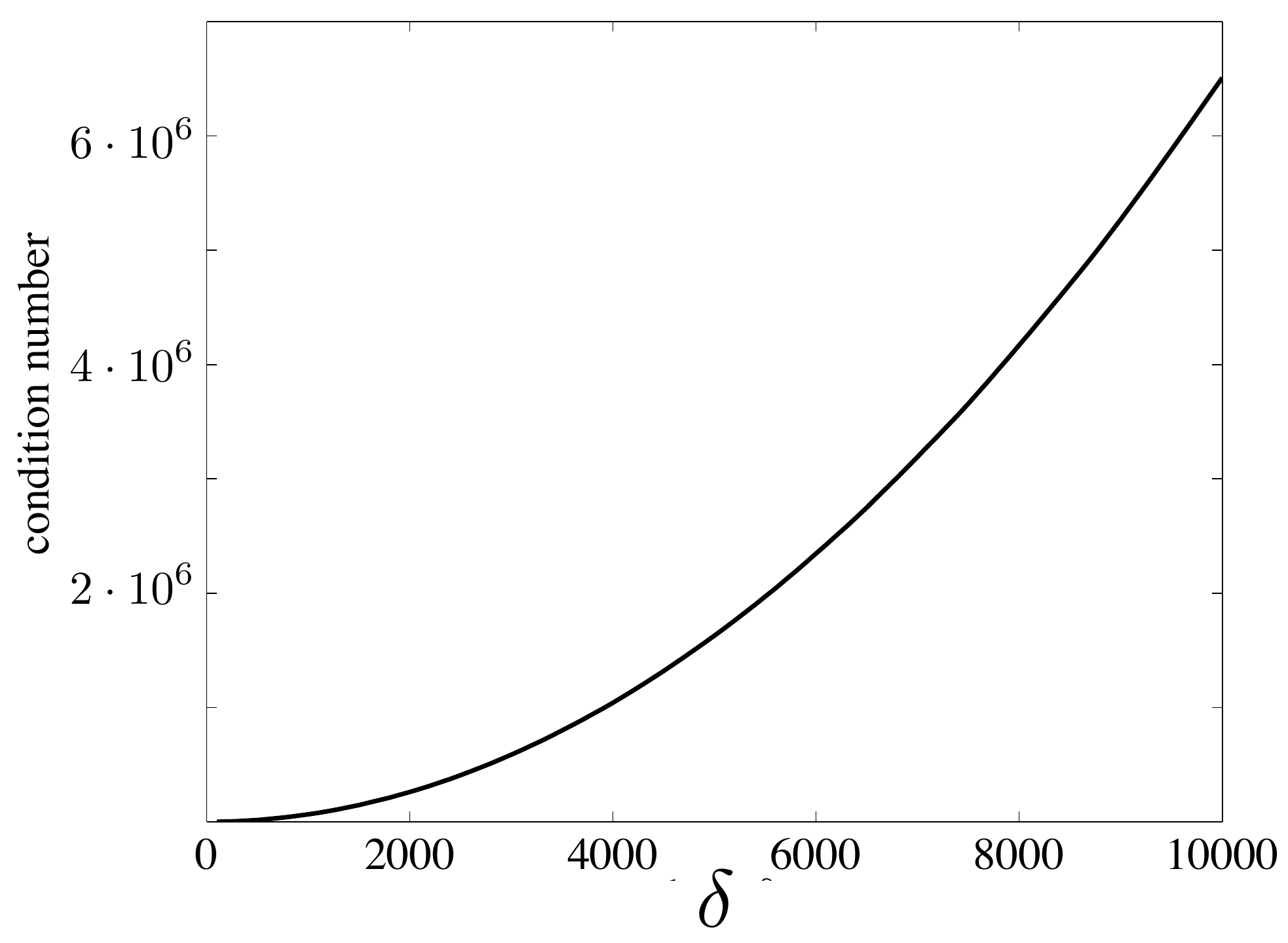} 
 \caption{$l^p$ condition numbers of $A_{1,p}(\epsilon_\delta)$ 
for $p=1$ (left), $p=2$ (center), and $p=\infty$ (right).
The slope of the internal layer is approximately $\delta$ and the 
thickness of the internal layer is approximately $1/\delta$.}
    \label{fig:op1conds}
\end{figure}

\begin{figure}[ht]
 \centering
    \includegraphics[width=.33\textwidth]{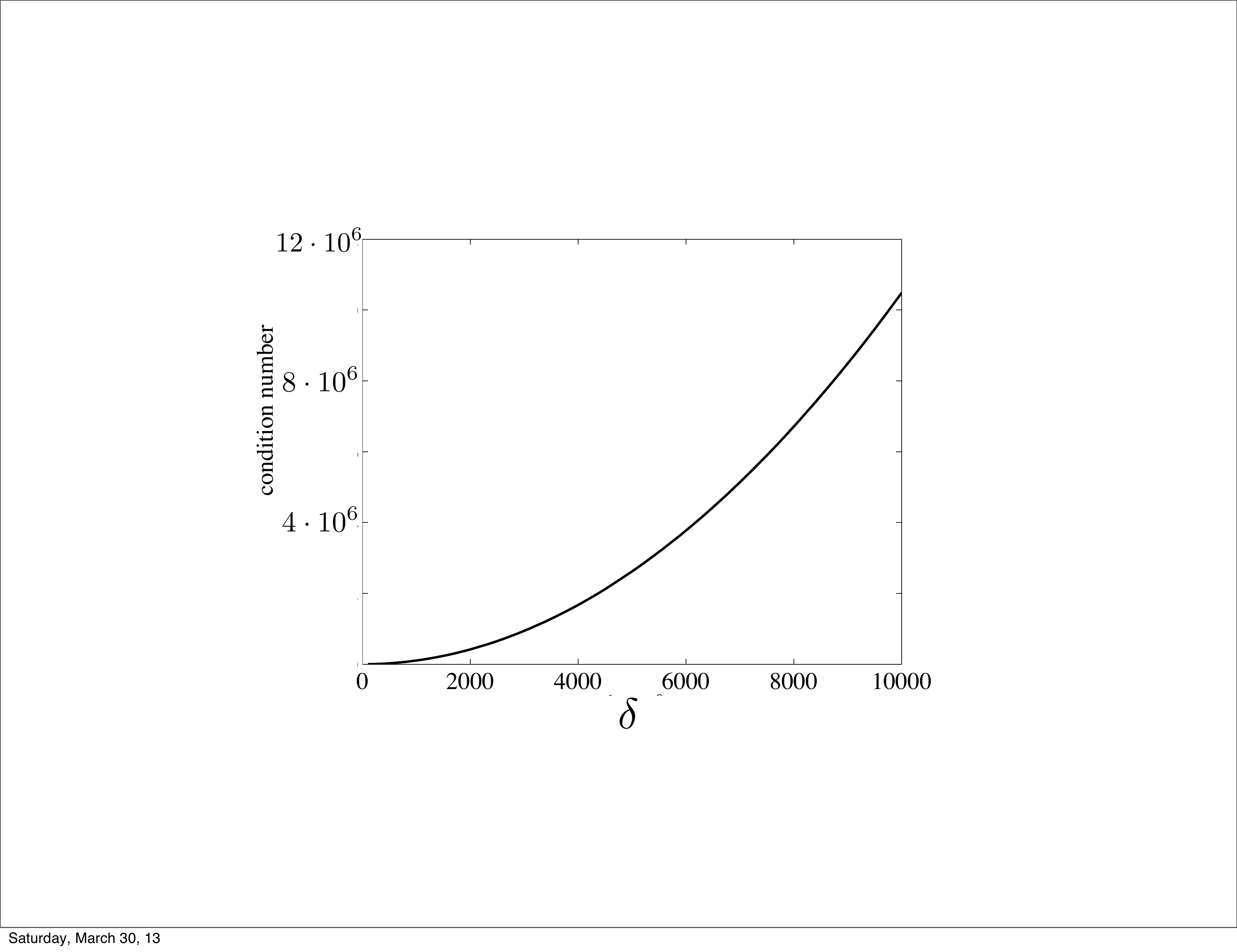}
    \includegraphics[width=.315\textwidth]{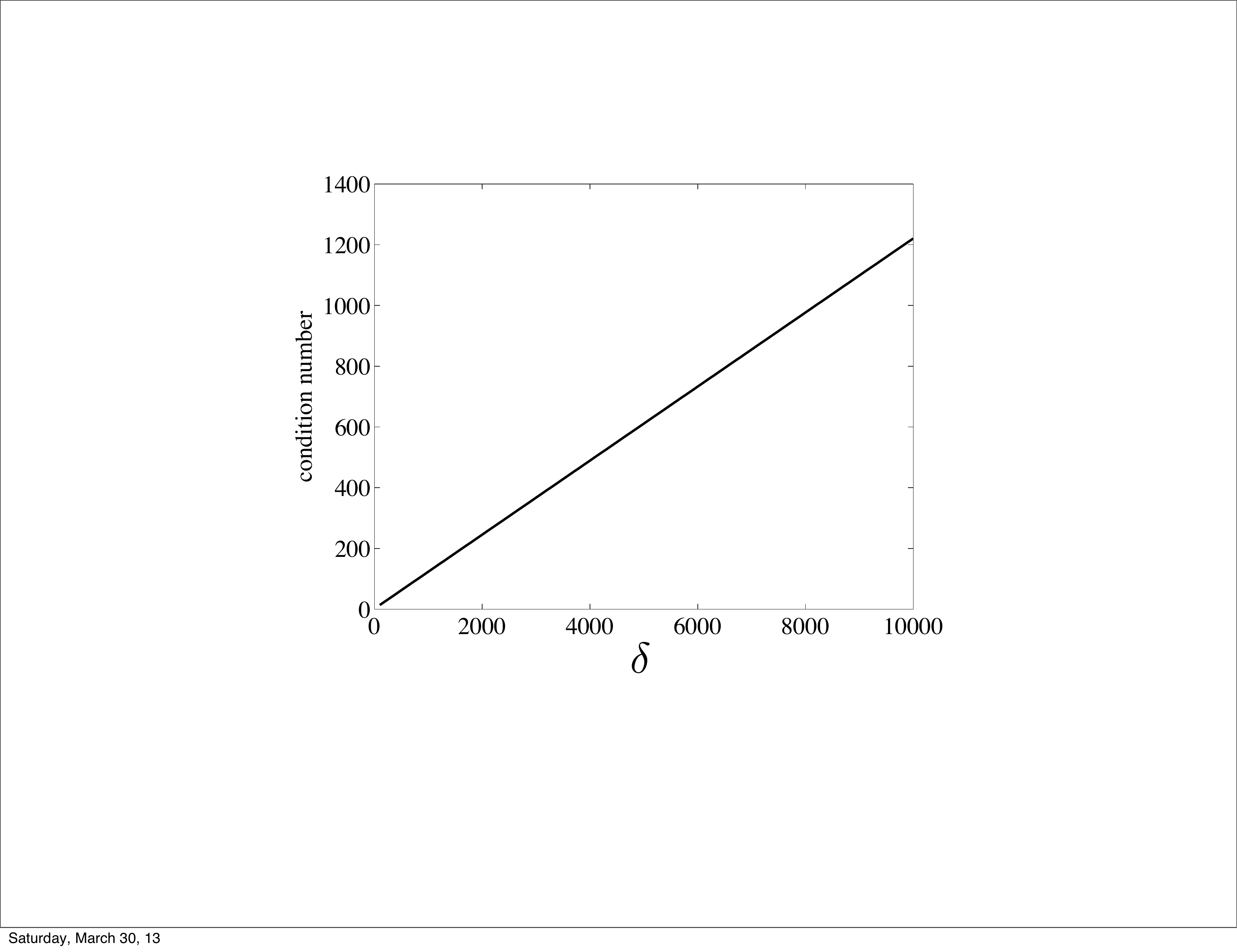}
    \includegraphics[width=.3\textwidth]{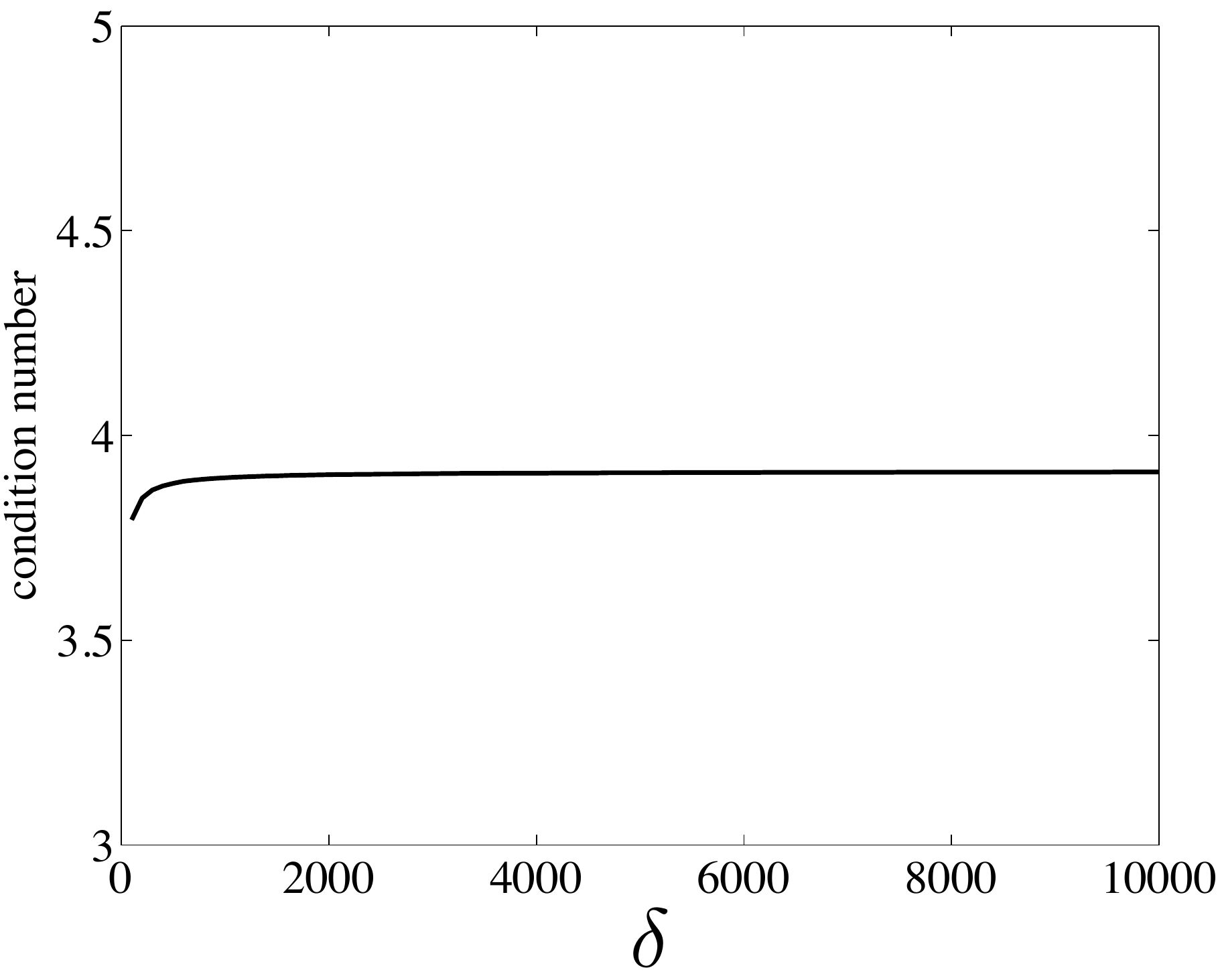} 
 \caption{$l^p$ condition numbers of $A_{2,p}(\epsilon_\delta)$ 
for $p=1$ (left), $p=2$ (center), and $p=\infty$ (right).}
    \label{fig:op2conds}
\end{figure}

In Figs. \ref{fig:op1conds} and \ref{fig:op2conds}, 
we plot numerical results for the family of 
functions $\epsilon_\delta$, where $\delta = 100j$, 
with $j=1,\ldots,100$. 
For each $\epsilon_\delta$, we formed the system matrices
for an adaptive norm-preserving discretization of the domain $[0,2]$ 
as described above. The 
$l^p$ condition numbers were computed by brute force 
(using the singular value decomposition in MATLAB).

We see  from the data that the condition numbers of the discrete operators 
do, indeed, exhibit the scaling properties expected from our analysis
of the continuous operators. Note that the 1-norm-preserving scheme to 
discretize (\ref{inteq1}) and the $\infty$-norm-preserving scheme to 
discretize (\ref{inteq2}) result in very well-conditioned matrices,
independent of the steepness of the internal layer.

\subsection{Convergence behavior using GMRES}

As discussed in section \ref{adaptivity}, it is reasonable to 
ask how standard iterative schemes work when applied to 
$l^p$-norm-preserving discretizations. We use GMRES here,
whose convergence behavior depends formally on the $l^2$ 
condition number of the 
system matrix. It is reasonable to  expect that the better conditioned 
systems (the 1-norm-preserving system for $A_1(\epsilon)$ and 
the $\infty$-norm-preserving system for $A_2(\epsilon)$) will 
fare better.

\begin{figure}[ht]
 \centering
  \includegraphics[width=.45\textwidth]{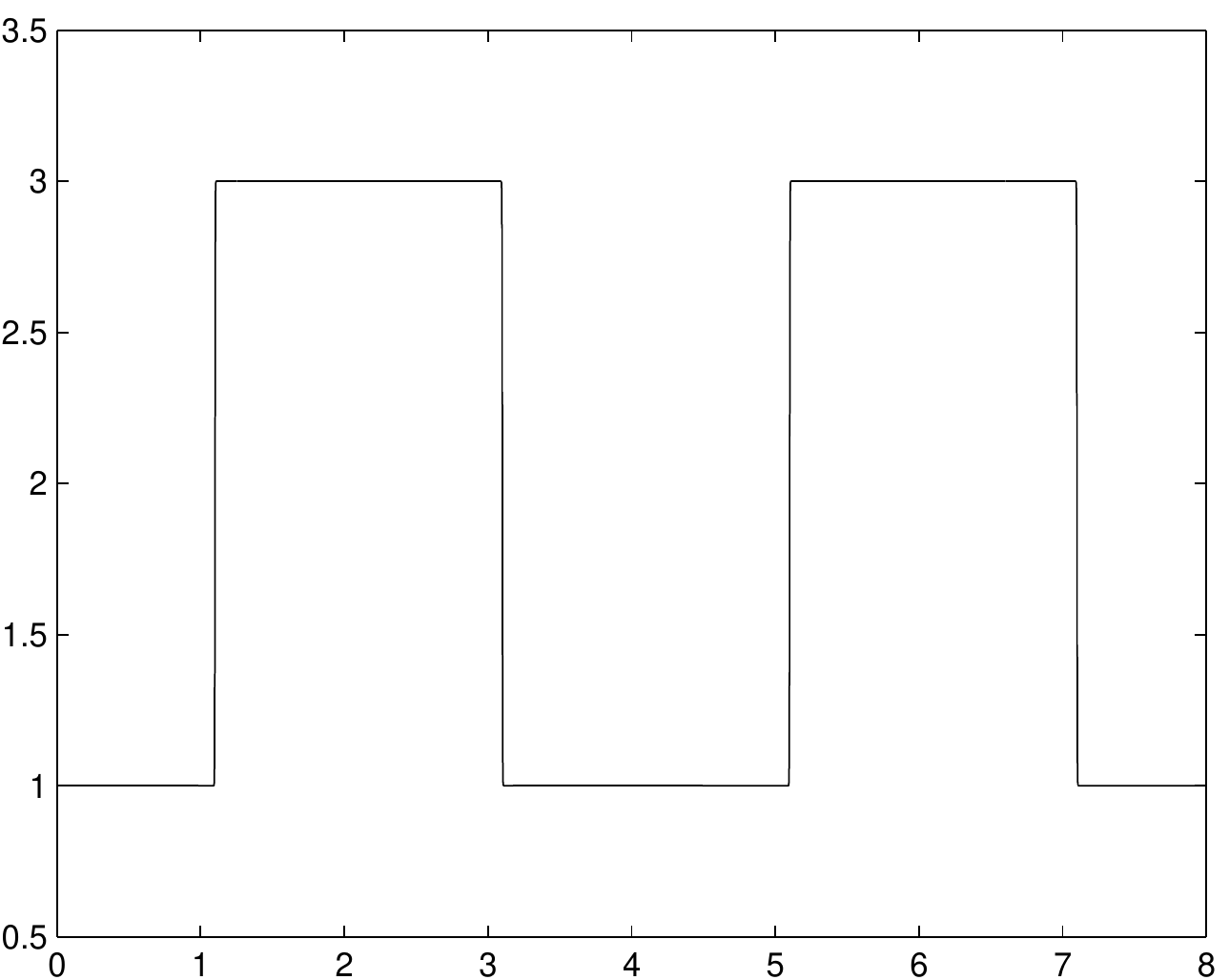}
  \includegraphics[width=.45\textwidth]{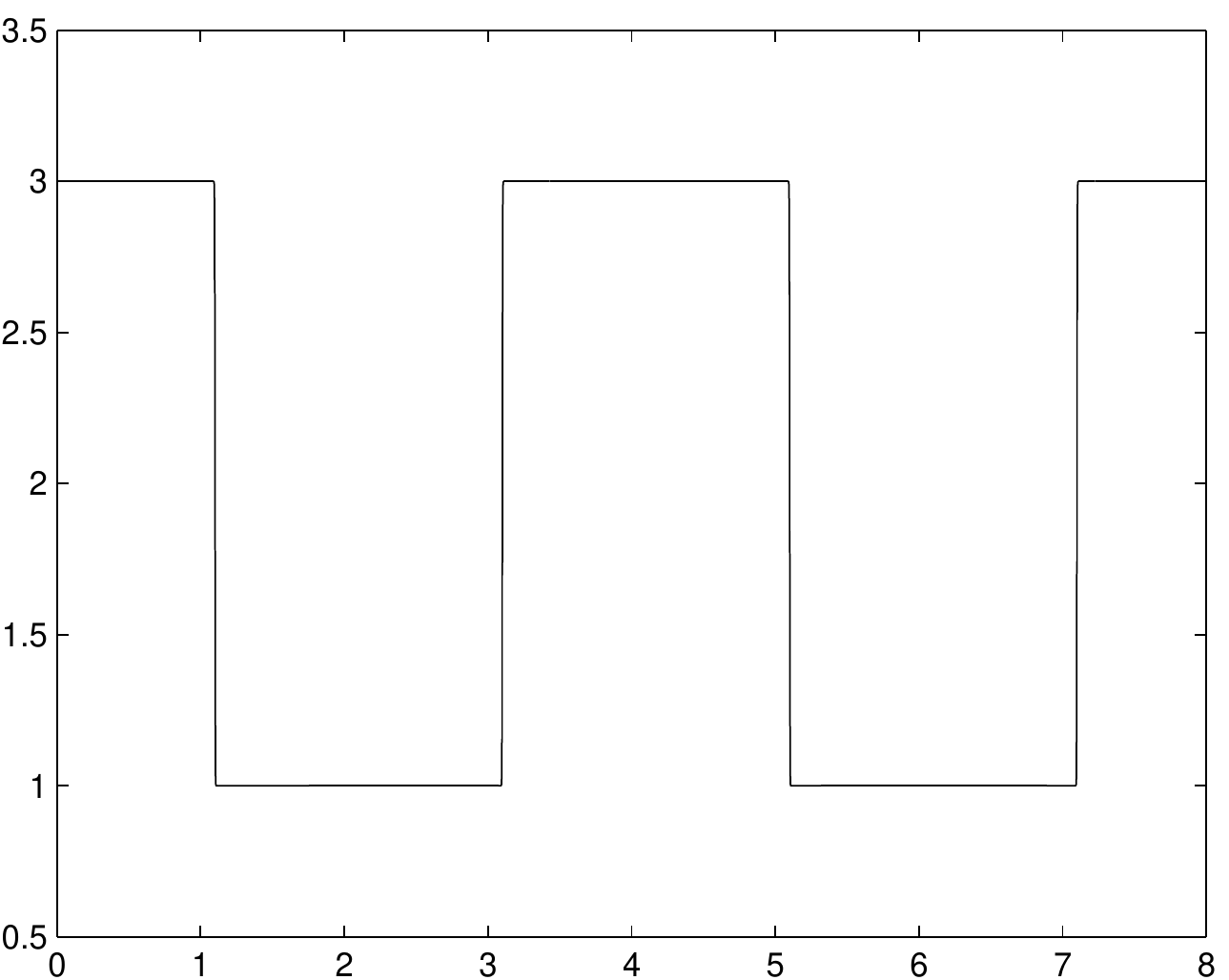}
 \caption{The ``double hill'' (left) and ``double well'' (right) functions
$\epsilon(x)$.}
 \label{fig:epsplots}
\end{figure}

For these experiments, we solve the ODE (\ref{ode1}), i.e.
\begin{equation}
 \frac{\partial}{\partial x} \left ( \epsilon (x) \frac{\partial u}{\partial x} \right ) = f\nonumber
\end{equation}
subject to inhomogeneous Dirichlet conditions, $u(a) = \gamma_a$ and 
$u(b) = \gamma_b$. If we let $l(x) = mx+c$ be a linear function satisfying 
the boundary conditions, then $v = u-l$ satisfies homogeneous
Dirichlet conditions and the ODE with a modificed right-hand side:
\[
 \frac{\partial}{\partial x} \left ( \epsilon (x) 
\frac{\partial v}{\partial x} \right ) = f - m \epsilon_x. 
\]
This problem can be addressed using one of the integral equations 
(\ref{inteq1}) or (\ref{inteq2}), from which the solution to the 
original problem is $u = v+l$. 
Here, we consider $f \equiv 1$, $\gamma_a = 1$ and $\gamma_b=2$. 
We consider two types of functions $\epsilon(x)$ 
that contain multiple internal layers by adding together 
several hyperbolic tangent functions, 
as in (\ref{epseq}), with multiple centers and $\delta=500$,
as shown in Fig. \ref{fig:epsplots}.
We refer to the left-hand profile as a ``double hill'' and the right-hand
profile as a ``double well''. 

\begin{figure}[ht]
\centering
\includegraphics[width=.45\textwidth]{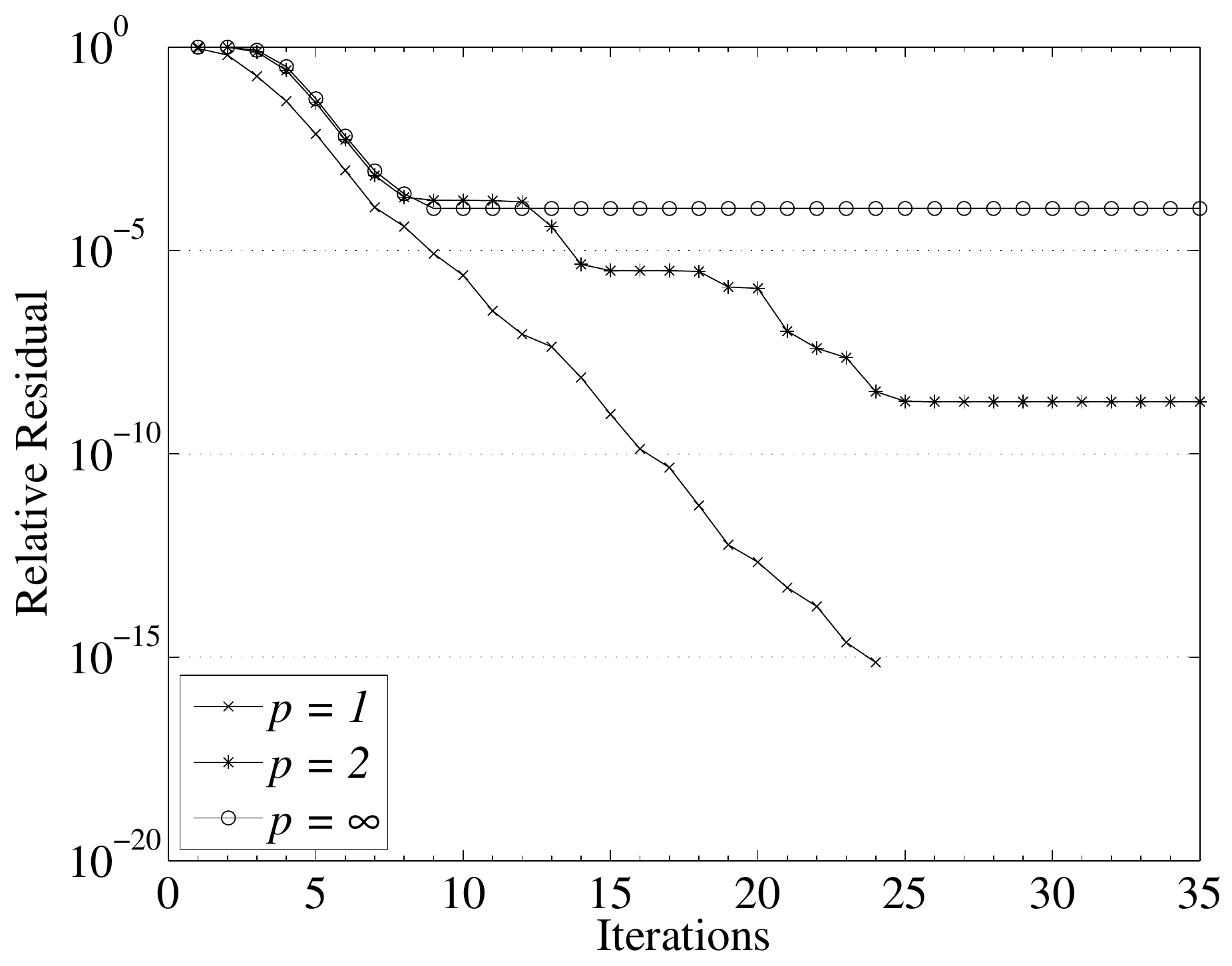}
\includegraphics[width=.45\textwidth]{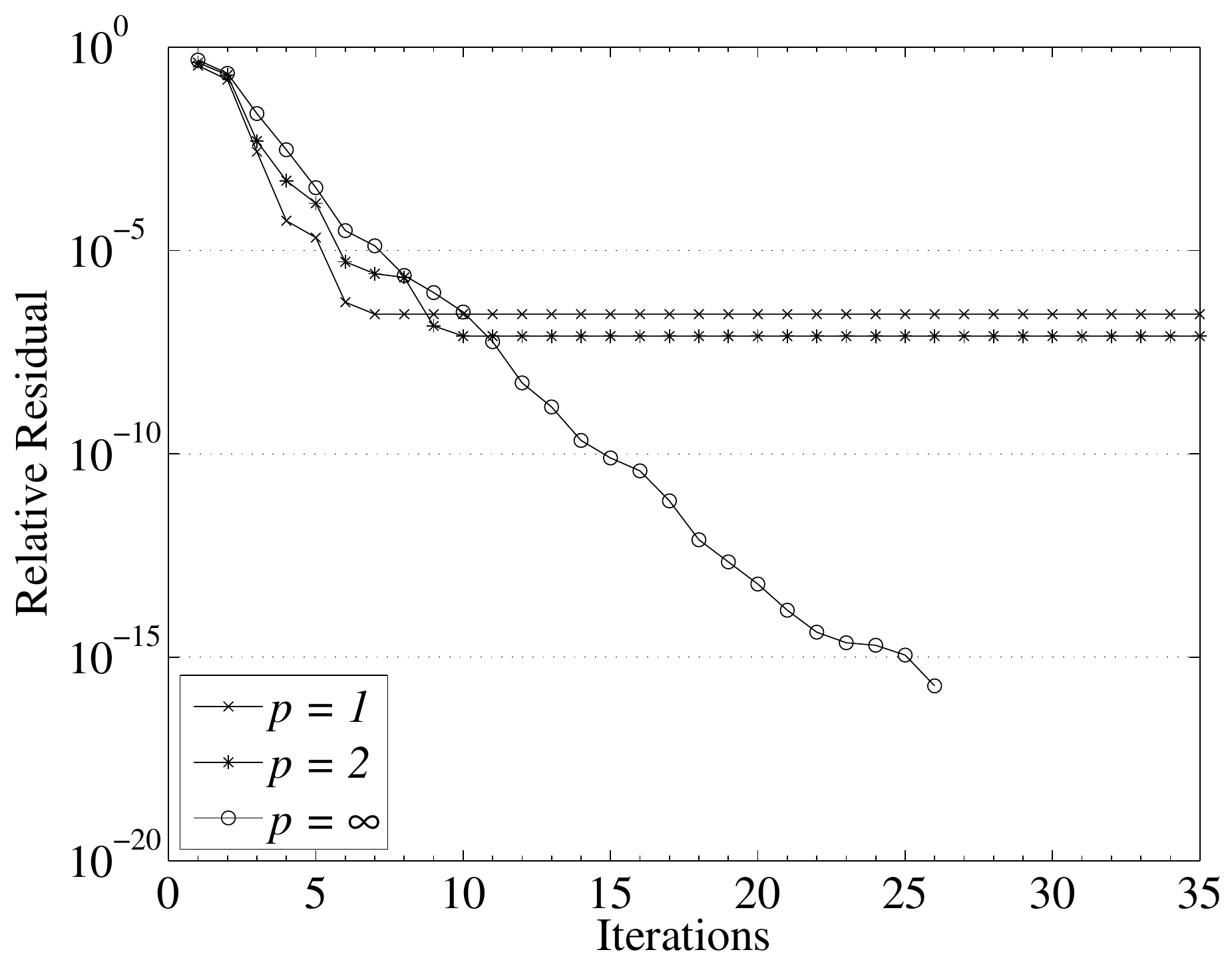}
\caption{Convergence of GMRES for the ``double hill'' $\epsilon(x)$.
The relative residual of the error at each iteration  is shown
using the $A_{1,p}(\epsilon)$ operator (left) and the $A_{2,p}(\epsilon)$ 
operator (right).}
\label{fig:coeff1}
 \end{figure}

\begin{figure}[ht]
\centering
\includegraphics[width=.45\textwidth]{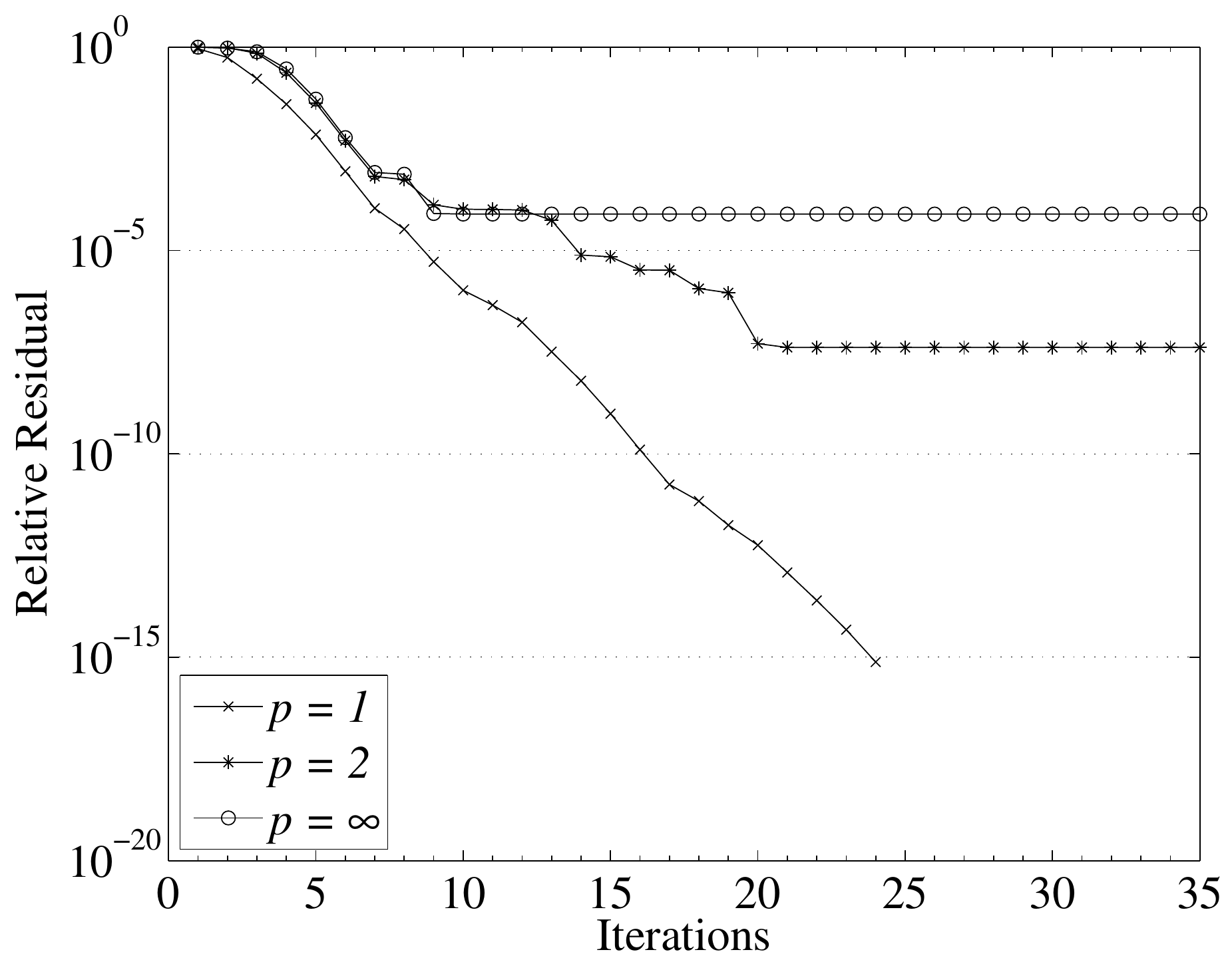}
\includegraphics[width=.45\textwidth]{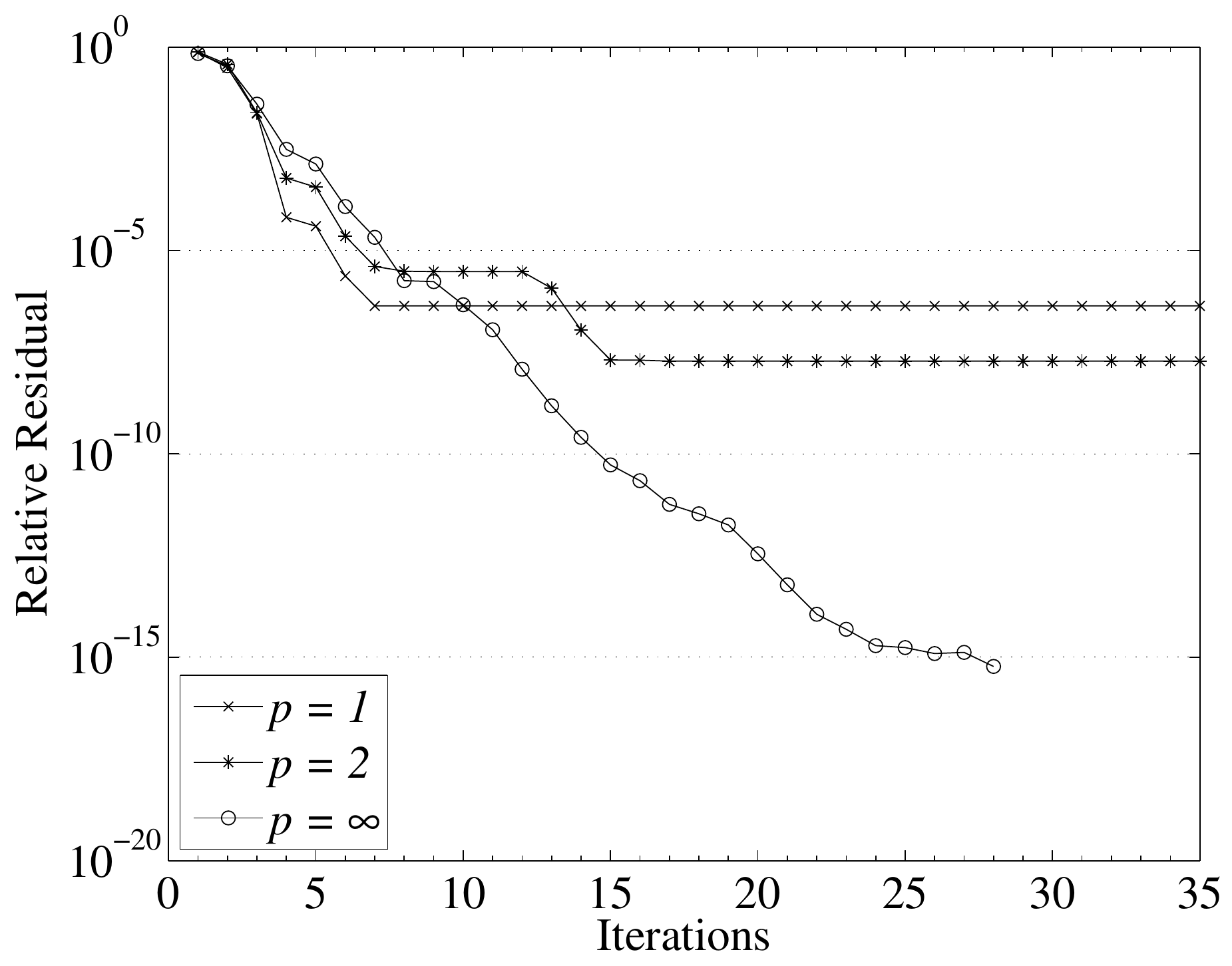}
\caption{Convergence of GMRES for the ``double well'' $\epsilon(x)$. 
The relative residual of the error at each iteration is shown
using the $A_{1,p}(\epsilon)$ operator (left) and the $A_{2,p}(\epsilon)$ 
operator (right).}
\label{fig:coeff2}
 \end{figure}

Using adaptive refinement, we obtain linear systems (\ref{linsystem1}) 
and (\ref{linsystem2}) as described above, for $p=1$, 2, and $\infty$. 
We solve the systems using GMRES and record the relative residuals for 
each step in Figs. \ref{fig:coeff1} and \ref{fig:coeff2}. 
The {\em $l^2$ condition numbers} of the discrete operators are shown in 
Table \ref{tab:conds}.

\begin{table}[ht]
\centering
 \begin{tabular}{|c|r|r|r|r|r|r|}
  \hline
  $\epsilon(x)$ & $A_{1,1}(\epsilon)$ & $A_{1,2}(\epsilon)$ & $A_{1,\infty} (\epsilon)$ & $A_{2,1}(\epsilon)$ & $A_{2,2} (\epsilon)$ & $A_{2,\infty}(\epsilon)$ \\ \hline
  ``Double Hill'' & 35.1453 & 979.052  &  86459.5 & 116010 & 978.240 & 31.1643 \\ \hline
  ``Double Well'' & 33.1648 & 977.744 & 98620.1 & 147328 & 977.411 & 27.9858 \\ \hline
 \end{tabular}
\caption{$l^2$ condition numbers for the discretized $A_{1,p}(\epsilon)$ and 
$A_{2,p}(\epsilon)$ operators.
 \label{tab:conds}}
\end{table}

Note that the $l^2$ condition numbers for $A_{1,1}(\epsilon)$ and 
$A_{2,\infty}(\epsilon)$ operators are the smallest, as expected.
Note also that these linear systems are solved much more easily
using GMRES. The other discretizations fail to reach the desired tolerance 
($10^{-15}$) in a reasonable number of iterations.
 
 \section{Discussion} \label{sec:discussion}

Our work in this paper was motivated by the observation that 
boundary integral equations are extremely robust when solving problems
of the type (\ref{pde1}) when $\epsilon$ is piecewise constant. 
In particular, a {\em charge} distribution on the dielectric
interface leads to 
well-conditioned integral equations involving the single layer potential
\cite{glee,gmoura,helsing2,martinsson,white}).
That charge density, however, is not a smooth function in the ambient space - 
it is a singular function supported on the interface alone.

In the variable
coefficient case, setting the unknown to be 
$\sigma = \Delta u$, as
in (\ref{inteq1}), corresponds to seeking the solution
in terms of a {\em volume} charge distribution. As
the internal layer becomes steeper and steeper, 
the function $\sigma(x)$ blows up, since it is converging to a
distribution and not a bounded function.
One interpretation of the $L^1$ norm-preserving discretization 
is that, in the discontinuous limit, the 
$l^1$-scaled unknown 
approximates the {\em strength} of the $\delta$-function along the
steep interface, rather than trying
to sample the $\delta$-function itself.

One concern with using the integral equation (\ref{inteq1})
is that we are only guaranteed tight bounds on accuracy in $L^1$,
using the standard estimate 
\[
\frac{\| e \|_1}{\| x \|_1} \leq \mbox{cond}_1(A_1)
 \frac{\|r \|_1}{\| b \|_1} \, 
\]
where $\tilde{x}$ is an approximate solution, $e = x - \tilde{x}$, 
and $r = A_1 \tilde{x} - b$ is the residual.
(This estimate applies to invertible Fredholm equations 
of the second kind as well as to finite-dimensional linear systems). 
Fortunately, the quantities of interest $u,u_x$ are computed as 
integral functionals of $\sigma$ using the representation
(\ref{urepsigma}) and are obtained with high accuracy.
The integral equation (\ref{inteq2}) can be discretized naively, 
corresponding, as noted earlier, to norm-preservation in $l^\infty$. 
While in some respects simpler, derivative data ($u_x$) must then 
be computed numerically. 

We are currently working on the extension of our analysis to 
higher-dimensional problems, and will report on the performance 
of such solvers at a later date. 

\section*{Acknowledgements}
This work was supported in part by the Applied
Mathematical Sciences Program of the U.S. Department of Energy
under Contract DEFGO288ER25053 and in part
by the Air Force Office of Scientific Research under 
NSSEFF Program Award FA9550-10-1-0180.

\appendix 
\section{Proof of Theorem~\ref{condthm1}}

Let $\mathcal{E}$ be a family of functions 
satisfying Properties 1 and 2
from Definition \ref{propdef}.
Let $\epsilon$ be an arbitrary function in $\mathcal{E}$ 
and let $A_1$ be given by
\[ A_1 \sigma(x) = (I+K_1) \sigma(x) = \sigma(x) + \frac{\epsilon_x(x)}{\epsilon(x)} \int G_x(x,y) \sigma(y) \, dy \, . \]
We now establish bounds for $A_1$ as on operator on 
$L^\infty [0,1] \cap C[0,1]$. To begin, we note that $|G_x(x,y)|$ is bounded by 1. Thus,
\begin{align}
\|A_1\|_\infty &= \sup_{ \|\sigma\|_\infty = 1} \sup_{x\in [0,1]} \left | 
\sigma(x) + \frac{\epsilon_x(x)}{\epsilon(x)} 
\int G_x(x,y) \sigma(y) \, dy \right | \\
& \leq 1 + \left \|\frac{\epsilon_x}{\epsilon} 
\right \|_\infty \, .
\end{align}
Let $x_*$ be the maximizer of $|\epsilon_x/\epsilon|$ and 
define the functions $\sigma_n$ by
\[ \sigma_n(y) = \left \{ \begin{array}{rcl} 1 & \mbox{ if } & y \leq x_* \\
                          1-2n(y-x_*) & \mbox{ if } & x_* < y < x_*+1/n \\
                          -1 & \mbox{ if } & y \geq x_* + 1/n
                          \end{array} \right. \, . 
\]
These functions are continuous and approximate the sign of 
$G_x(x_*,y)$. A straightforward computation shows that
\begin{align}
   \|A_1\|_\infty &= \sup_{ \|\sigma\|_\infty = 1} \sup_{x\in [0,1]} \left | 
\sigma(x) + \frac{\epsilon_x(x)}{\epsilon(x)} 
\int G_x(x,y) \sigma(y) \, dy \right | \\
   &\geq \sup_{x\in [0,1]} \left | \sigma_n(x) + \frac{\epsilon_x(x)}{\epsilon(x)} 
\int G_x(x,y) \sigma_n(y) \, dy \right |\\
   &\geq \left| \frac{\epsilon_x(x_*)}{\epsilon(x_*)} \right| \left ( 
\int |G_x(x_*,y)| \, dy - \frac{2}{n} \right ) - \sigma_n(x_*) \\
   &\geq \left \| \frac{\epsilon_x}{\epsilon} \right \|_\infty 
\left ( \frac{1}{4} - \frac{2}{n} \right ) -1 
\end{align}
so that 
\[
\|A_1\|_\infty \geq \frac{1}{4} \left \| 
\frac{\epsilon_x}{\epsilon} \right \|_\infty - 1 \, .
\]

We note that $A_1^{-1}$ is given by

\[ A_1^{-1} g(x) = (I-R_1)g(x) = g(x) - \frac{\epsilon_x(x)}{\epsilon(x)^2} \left ( \int_0^x g(t) \epsilon(t) \, dt  - \frac{\int_0^1 \frac{1}{\epsilon(s)} \int_0^s g(t) \epsilon(t) \, dt \, ds} { \int_0^1 \frac{1}{\epsilon(s)} \, ds} \right ) \, .
\]
It is straightforward to see that
  
\begin{align}
  \|A_1^{-1}\|_\infty &= \sup_{ \|g\|_\infty = 1 } \| (I - R_1) g \|_\infty \\
&\leq 1 + \sup_{ \|g\|_\infty = 1 } \sup_{x\in [0,1]} \left |\frac{\epsilon_x(x)}{\epsilon(x)^2} 
\left( \int_0^x g(t) \epsilon(t) \, dt  - 
\frac{\int_0^1 \frac{1}{\epsilon(s)} \int_0^s g(t) \epsilon(t) \, dt \, ds}
{ \int_0^1 \frac{1}{\epsilon(s)} \, ds} \right ) \right | \\
   &\leq 1 + \left \|\frac{\epsilon_x}{\epsilon} \right\|_\infty  \frac{1}{m}
\left ( 1 + \frac{M}{m} \right ) \|\epsilon \|_1 \, .
  \end{align}
Again, let $x_*$ be the maximizerof $| \epsilon_x/ \epsilon |$ and let $m_\epsilon$ be the minimum of $\epsilon$ on $[0,1]$. We define the function $g_\epsilon$ as follows
\[ g_\epsilon(x) = \left \{ \begin{array}{rcl} m_\epsilon/\epsilon(x) & \mbox{ if } 
& x \leq x_*/2 \\
-m_\epsilon/\epsilon(x) & \mbox{ if } 
& x_*/2 < x \leq x_* \\
m_\epsilon/\epsilon(x) & \mbox{ if } 
& x_* < x \leq (1+x_*)/2 \\
-m_\epsilon/\epsilon(x) & \mbox{ if } 
& (1+x_*)/2 < x \leq 1 \end{array} \right. \, . 
\]
The function $g_\epsilon$ is such that the integral $\int_0^x g_\epsilon(t) \epsilon(t) \, dt$ is zero at $x=x^*,0,$ and $1$ and positive otherwise. Let $g_n$ be continuous functions which satisfy $\|g_n\|_\infty = 1$ and converge pointwise to $g_\epsilon$. A few straightforward computations and an application of the dominated convergence theorem yield
\begin{align}
 \|A_1^{-1}\|_\infty &= \sup_{ \|g\|_\infty = 1 } \left \| (I-R_1)g \right \|_\infty \\
    &\geq \lim_{n \to \infty}  \left \| (I-R_1)g_n \right \|_\infty \\
    &\geq \lim_{n \to \infty}  \left | R_1 g_n (x_*) \right | - 1\\
    &\geq \left | \frac{\epsilon_x(x_*)}{\epsilon(x_*)^2}\right | 
\lim_{n\to \infty} \frac{\int_0^1 \frac{1}{\epsilon(s)} 
\int_0^s g_n(t) \epsilon(t) \, dt \, ds}
{ \int_0^1 \frac{1}{\epsilon(s)} \, ds}  - \left | 
\frac{\epsilon_x(x_*)}{\epsilon(x_*)^2} \right| \lim_{n\to \infty}  
\int_0^{x_*} g_n(t) \epsilon(t) \, dt - 1 \\
    &= \left | \frac{\epsilon_x(x_*)}{\epsilon(x_*)^2}\right|
\frac{\int_{0}^1 \frac{1}{\epsilon(s)} \int_{0}^s g_\epsilon(t) \epsilon(t) \, dt \, ds}
{ \int_0^1 \frac{1}{\epsilon(s)} \, ds} - 1 \\
    &\geq \frac{1}{M} \left \| \frac{\epsilon_x}{\epsilon} 
\right \|_\infty \frac{m_\epsilon m}{8M} - 1 \\
    &\geq \left \| \frac{\epsilon_x}{\epsilon} \right \|_\infty 
\frac{m^2}{8M^2} - 1 \, .
\end{align}
  
We next establish bounds on $A_1$ as an operator on $L^p[0,1]$,
for $1 < p < \infty$. 

\begin{align}
   \|A_1\|_p &= \sup_{ \|\sigma\|_p = 1} \left \| (I+K_1) \sigma  \right \|_p \\
   &\leq  1 +  \sup_{ \|\sigma\|_p = 1} \left ( \int \left | 
\frac{\epsilon_x(x)}{\epsilon(x)}  \int G_x(x,y) \sigma(y) \, dy 
\right|^p \, dx \right )^{1/p} \\
   &\leq 1 +  \sup_{ \|\sigma\|_p = 1} \left \| \left(
\frac{\epsilon_x(x)}{\epsilon(x)} \right )^p \right \|_1^{1/p} \left\| 
\int \left |G_x(\cdot,y) \sigma(y) \right |^p \, dy \right \|_\infty^{1/p} \\
   &\leq 1 +  \left \|\frac{\epsilon_x(\cdot)}{\epsilon(\cdot)} \right\|_p \, .
\end{align}  
Because $\mathcal{E}$ satisfies Property 1, we may choose
$0 \leq \delta \ll 1$ and a neighborhood $V = B(\xi,c) \subset [0,1]$ 
centered at $\xi$ and of radius $c$ such that 
\[ \left \| \frac{\epsilon_x}{\epsilon} \cdot 1_V \right \|_p \leq \delta \left \| \frac{\epsilon_x}{\epsilon} \right \|_p \, .
\]
For concreteness, assume $c \leq 1/2$. 
Note now that a density $\sigma$ is the second derivative of a function 
$u$ with homogeneous Dirichlet boundary values. 
In particular, $u'(x) = \int G_x(x,y) \sigma(y) \, dy$. As a result, 
the function $u'$ integrates to zero (since $u(1) = u(0)=0$). 
This observation permits us to build densities $\sigma$ with the 
desired properties. In particular, we'd like a density $\sigma$ such 
that $u'(x) = \int G_x(x,y) \sigma(y) \, dy$ is small only in the 
neighborhood where $\epsilon_x/\epsilon$ satisfies the above property. 
We choose $u'$ to be of the form
\[  u'(x) = \left \{ \begin{array}{rcl}
                       a(x-\xi)^2 & \mbox{ for } & x \leq \xi\\
                       b(x-\xi)^2 & \mbox{ for } & x > \xi
                      \end{array} \right. \,  \]
where $a$ and $b$ are chosen such that $u'$ integrates to zero. 
Setting $b=1$ and $a = (\xi - 1)^3/ \xi^3$ is sufficient. 
This yields
\[ \sigma(x) = \left \{ \begin{array}{rcl}
               \frac{2 (\xi-1)^3}{\xi^3} (x-\xi) & \mbox{ for } & x \leq \xi \\
                    2 (x-\xi) & \mbox{ for } & x > \xi
                        \end{array} \right. \]
The $L_p$ norm of the above function satisfies
\[ \|\sigma\|_p \leq \|\sigma\|_\infty \leq \frac{2(1-c)}{c^3} \, . \]
Let $\sigma_\epsilon = \sigma / \|\sigma\|_p$. 
Then the corresponding $u'(x) = \int G_x(x,y)\sigma_\epsilon(y) \, dy$ 
is given by
\[
\int G_x(x,y) \sigma_\epsilon(y) \, dy = \left\{ 
       \begin{array}{rcl}
    \frac{1}{\|\sigma \|_p} \frac{(\xi-1)^3}{\xi^3} (x-\xi)^2 & \mbox{ for } 
             & x \leq \xi \\
    \frac{1}{\|\sigma \|_p} (x-\xi)^2 & \mbox{ for } & x > \xi
       \end{array} \right.  \, .
\]
This provides a minimum value of 
$|\int G_x(x,y) \sigma_\epsilon(y) \, dy|$ on $[0,1] \setminus V$ 
which satisfies
\[ \min_{x\in [0,1]\setminus V} \left | \int G_x(x,y) 
\sigma_\epsilon(y) \, dy \right | \geq \frac{c^8}{(1-c)^4} \, .
\]
We then have
\begin{align}
 \|A_1\|_p &= \sup_{ \|\sigma\|_p = 1} \left \| (I+K_1) \sigma \right \|_p \\
   &\geq  \sup_{ \|\sigma\|_p = 1} \left ( \int \left ( \frac{\epsilon_x(x)}{\epsilon(x)}  \int G_x(x,y) \sigma(y) \, dy \right )^p \, dx \right )^{1/p} - 1 \\
   &\geq \left ( \int \left ( \frac{\epsilon_x(x)}{\epsilon(x)}  \int G_x(x,y) \sigma_\epsilon(y) \, dy \right )^p \, dx \right )^{1/p} - 1 \\
   &\geq  (1-\delta) \frac{c^8}{1-c^4} \left \| \frac{\epsilon_x(\cdot)}{\epsilon(\cdot)} \right \|_p - 1 \\
   &= C(c,\delta) \left \| \frac{\epsilon_x(\cdot)}{\epsilon(\cdot)} \right \|_p - 1 \, .
\end{align}
Let $1 < p < \infty$ and $1/p + 1/q = 1$. Then
\begin{align}
   \|A_1^{-1}\|_p &=  \sup_{ \|g\|_p = 1 } \left \| (I-R_1) g \right \|_p \\
   &\leq 1 + \sup_{ \|g\|_p = 1} \|R_1 g\|_p \\
   &\leq 1 + \sup_{ \|g\|_p = 1 } \left (\int_0^1 \left | \frac{\epsilon_x(x)}{\epsilon(x)^2} \left ( \int_0^x g(t) \epsilon(t) \, dt  - \frac{\int_0^1 \frac{1}{\epsilon(s)} \int_0^s g(t) \epsilon(t) \, dt \, ds} { \int_0^1 \frac{1}{\epsilon(s)} \, ds} \right ) \right |^p \, dx \right )^{1/p}\\
   &\leq 1 +  \left \| \left ( \frac{\epsilon_x}{\epsilon^2} \right )^p \right \|_1^{1/p} \sup_{ \|g\|_p = 1 } \left ( \sup_{x \in [0,1]} \left | \int_0^x g(t) \epsilon(t) \, dt  - \frac{\int_0^1 \frac{1}{\epsilon(s)} \int_0^s g(t) \epsilon(t) \, dt \, ds} { \int_0^1 \frac{1}{\epsilon(s)} \, ds} \right |^p \right )^{1/p} \\
   &\leq 1 +\left \| \frac{\epsilon_x}{\epsilon^2} \right \|_p   \sup_{ \|g\|_p = 1 } \left (\int_0^1 |g(t) \epsilon(t)| \, dt  + \frac{\int_0^1 \frac{1}{\epsilon(s)} \int_0^1 |g(t) \epsilon(t)| \, dt \, ds} { \int_0^1 \frac{1}{\epsilon(s)} \, ds} \right ) \\
   &\leq 1 +  \left \| \frac{\epsilon_x}{\epsilon^2} \right \|_p \left ( 1 + \frac{M}{m} \right ) \sup_{ \|g\|_p = 1 } \| g \epsilon \|_1 \\
   &\leq 1 + \left \| \frac{\epsilon_x}{\epsilon^2} \right \|_p \left ( 1 + \frac{M}{m} \right ) \| \epsilon \|_q \\
   &\leq 1 +  \left \| \frac{\epsilon_x}{\epsilon} \right \|_p \left ( 1 + \frac{M}{m} \right ) \frac{M}{m} \, .
  \end{align} 

Let $V = B(\xi,c) \subset [0,1]$ as above. 
We define a function $g_\epsilon$ as follows:
\[ g_\epsilon(x) = \left \{ \begin{array}{rcl} 0 & \mbox{ if } & x \leq \xi - c \\
                              \frac{1}{\epsilon(x)} & \mbox{ if } & \xi-c < x \leq \xi \\
                              -\frac{1}{\epsilon(x)} & \mbox{ if } & \xi < x \leq \xi + c \\
                              0 & \mbox{ if } & x > \xi+c 
                             \end{array} \right. \, .
\]
It is easy to see that
\[ \int_0^x g_\epsilon(t) \epsilon(t) \, dt = \left\{ 
  \begin{array}{rcl} 0 & \mbox{ if } & x \leq \xi - c \\
                     x- \xi + c & \mbox{ if } & \xi-c < x \leq \xi \\
                     \xi-x+c & \mbox{ if } & \xi < x \leq \xi + c \\
                      0 & \mbox{ if } & x > \xi+c 
  \end{array} \right. \, ,
\]
that
\[ \int_0^1 \frac{1}{\epsilon(s)} \int_0^s g(t) \epsilon(t) \, dt \, ds 
\geq \frac{2}{M} \int_0^c t \, dt = \frac{c^2}{M}  \, ,
\]  
and that
\[  \|g_\epsilon \|_p \leq \frac{(2c)^{1/p}}{m} \leq \frac{1}{m} \, . \]
>From these,
\begin{align}
    \|A_1^{-1}\|_p &=  \sup_{ \|g\|_p = 1 } \left \| (I-R_1)g \right \|_p \\
    &\geq \sup_{ \|g\|_p = 1 } \left \| R_1 g \right \|_p - 1 \\
    &\geq \frac{1}{ \|g_\epsilon\|_p} \left ( \int_0^1 \left | \frac{\epsilon_x(x)}{\epsilon(x)^2} \left ( \int_0^x g_\epsilon(t) \epsilon(t) \, dt  - \frac{\int_0^1 \frac{1}{\epsilon(s)} \int_0^s g_\epsilon(t) \epsilon(t) \, dt \, ds} { \int_0^1 \frac{1}{\epsilon(s)} \, ds} \right ) \right |^p \, dx \right) ^{1/p} - 1 \\
    &\geq \frac{1}{ \|g_\epsilon\|_p} \left ( \int_{[0,1] \setminus V} \left | \frac{\epsilon_x(x)}{\epsilon(x)^2} \left ( \int_0^x g_\epsilon(t) \epsilon(t) \, dt  - \frac{\int_0^1 \frac{1}{\epsilon(s)} \int_0^s g_\epsilon(t) \epsilon(t) \, dt \, ds} { \int_0^1 \frac{1}{\epsilon(s)} \, ds} \right ) \right |^p \, dx \right) ^{1/p} - 1 \\
    &= \frac{1}{ \|g_\epsilon\|_p} \left ( \int_{[0,1] \setminus V} \left | \frac{\epsilon_x(x)}{\epsilon(x)^2} \left ( \frac{\int_0^1 \frac{1}{\epsilon(s)} \int_0^s g_\epsilon(t) \epsilon(t) \, dt \, ds} { \int_0^1 \frac{1}{\epsilon(s)} \, ds} \right ) \right |^p \, dx \right) ^{1/p} - 1\\
    &\geq  \left \| \frac{\epsilon_x}{\epsilon^2}\right \|_p m (1-\delta) \frac{c^2}{M}  - 1 \\
    &\geq \left \| \frac{\epsilon_x}{\epsilon}\right \|_p (1-\delta) \frac{mc^2}{M^2}  - 1 \, .
  \end{align}
>From the above, we see that there exist constants $C'_1$ and $C'_2$ -- depending only on $m,M,c,$ and $\delta$ -- such that 

\begin{align*}
 C_1'  \left \| \frac{\epsilon_x}{\epsilon}\right \|_p - 1 &\leq \|A_1\|_p \leq C_2'  \left \| \frac{\epsilon_x}{\epsilon}\right \|_p + 1 \\
 C_1'  \left \| \frac{\epsilon_x}{\epsilon}\right \|_p - 1 &\leq \|A^{-1}_1\|_p \leq C_2'  \left \| \frac{\epsilon_x}{\epsilon}\right \|_p + 1 \, ,
\end{align*}
so that there are constants $C_1$ and $C_2$ -- depending only on $m,M,c,$ and $\delta$ -- such that

\[  C_1  \left \| \frac{\epsilon_x}{\epsilon}\right \|_p - 1 \leq \mbox{cond}_p(A_1) \leq C_2  \left \| \frac{\epsilon_x}{\epsilon}\right \|_p + 1 \, , \]
which completes the proof. \qedhere
     
\bibliographystyle{plain}

\end{document}